
\input amstex

\documentstyle{amsppt}

\magnification=\magstep 1
\hsize=15,7truecm
\vsize=22truecm
\hoffset=0,3truecm
\voffset=1truecm

\topmatter

\title 
The Bianchi--Darboux transform of $L$-isothermic surfaces
\endtitle

\author  Emilio Musso  and  Lorenzo Nicolodi \endauthor

\abstract
We study an analogue of the classical B\"acklund transformation 
for $L$-i\-so\-thermic surfaces in Laguerre geometry, the Bianchi--Darboux 
transformation. We show how to construct the Bianchi--Darboux transforms
of an $L$-isothermic surface 
by solving an integrable linear 
differential system. We then establish a permutability theorem for iterated 
Bianchi--Darboux transforms. 
\endabstract

\address  Dipartimento di Matematica Pura ed Applicata,
Universit\`a di L'Aquila, via Vetoio, I-67010 Coppito (L'Aquila),
Italy 
\endaddress

\email musso\@univaq.it \endemail 

\address   Dipartimento di Matematica "G. Castelnuovo",
Universit\`a di Roma "La Sapienza", piazzale A. Moro 2, I-00185 Roma,
Italy
\endaddress 

\email nicolodi\@mat.uniroma1.it\endemail


\thanks Partially supported by MURST Project ``Propriet\`a geometriche delle
variet\`a reali e com\-ples\-se''.
\endthanks  



\endtopmatter

\document

\subheading{1. Introduction}  

\bigskip

\flushpar Certain types of integrable
non-linear PDEs (soliton equations) arise in differential geometry
as compatibility conditions for the linear equations obeyed by 
frames adapted to surfaces in higher dimensional manifolds. 
In a number of situations, the construction of new solutions of the arising
PDE relies on the existence of B\"acklund type transformations 
for the surfaces and on their permutability properties. 
Well-known examples include
pseudo-spherical 
surfaces, affine minimal surfaces, and isothermic surfaces in 
M\"obius geometry 
[CT],[BHPP],[CGS],[TT].
The loop group approach to soliton 
theory, recently developed by Terng--Uhlenbeck and others [TU], explains
the uniformity of results in several of these examples.

\medskip
\flushpar This paper concentrates on $L$-isothermic surfaces. 
An immersion $f : M \to \Bbb R^3$ is called $L$-isothermic if about each point
of $M$ there exist curvature line coordinates $(x,y)$ which are conformal 
with respect to
the third fundamental form of $f$. 
$L$-isothermic surfaces are invariant under a group of contact 
transformations, the Laguerre group, which is isomorphic to the 
identity component of the 
isometry group of Minkowski 4-space.
This group acts transitively on the set of 
oriented 2-spheres 
and points as well as on the set of
oriented planes of $\Bbb R^3$. 
We first became interested in $L$-isothermic surfaces when 
studying the 
(infinitesimal) deformation problem for surfaces under the Laguerre group.
Namely, $L$-isothermic surfaces are the only surfaces allowing 
1-parameter families of second order deformations [MN1],[MN2]. 
Such families correspond to the solutions of the
non-linear fourth-order equation $\Delta(w_{xy}w^{-1})=0$ which in turn is
equivalent to the Gauss-Codazzi equations
and arises
as the integrability 
condition of a linear differential system containing a free parameter [MN3].

\medskip
\flushpar The theory of $L$-isothermic surfaces presents many analogies  
with that of isothermic surfaces in M\"obius geometry. The latter has 
received much attention in recent studies [Bu],[BHPP],[He],[HP],[HMN], 
especially in relation with the general theory of curved flats in 
pseudo-Riemannian symmetric spaces [FP].  
Curved flats arise in 1-parameter families as solutions to a certain 
integrable system expressed in the form of a zero-curvature equation 
with spectral parameter.
This approach to the study of isothermic surfaces provides a uniform 
framework for understanding their theory of transformations and deformation.
The notion of isothermic submanifolds have also been defined in other
contexts; for instance, a general theory of isothermic submanifolds in 
symmetric $R$-spaces have been recently developed in [BPP].

\medskip
\flushpar In this paper, we study 
a geometric transformation for $L$-isothermic 
surfaces by realizing them as enveloping surfaces of a 2-sphere congruence. 
The Laguerre invariant conditions 
that the congruence preserve curvature lines and that the third fundamental 
forms of the two envelopes at corresponding points be conformal yield 
that both surfaces are $L$-isothermic. This is the content of Theorem 1.
The resulting
congruence is an analogue of the Darboux congruence occurring 
in M\"obius geometry and its two envelopes are said to form a 
Bianchi--Darboux pair. A Bianchi--Darboux pair gives rise to a 
curved flat in the Grassmannian $\tilde{G}_{1,1}(\Bbb R^4_1)$ of oriented 
2-planes of signature $(1,1)$ in 
$\Bbb R^4_1$; conversely, a curved flat in this Grassmannian 
only determines the normals of a Bianchi--Darboux pair. 
Moreover, the spectral
deformation of curved flats amounts to second order deformation as in the
conformal situation.
In Section 5, we 
prove the existence of Bianchi--Darboux transforms by explicitly constructing 
new $L$-isothermic surfaces from a given 
one; the construction requires solving an integrable 
linear differential system. This furnishes a geometrical method 
of deriving new solutions of the defining PDE 
$\Delta(w_{xy}w^{-1})=0$ from any given one. 
Finally, in Section 6, we establish a permutability theorem for the 
Bianchi--Darboux transformation.

\medskip
\remark{Acknowledgements} It is a pleasure to thank Fran Burstall and
Udo Hertrich-Jeromin for many interesting and instructive discussions 
concerning isothermic surfaces and their transformation theory.
Thanks also to the referee for his comments and questions.
\endremark

\bigskip

\subheading{2. $L$-isothermic surfaces}

\bigskip
 
\flushpar Let $f : M\to \Bbb R^3$ be an immersion of a surface $M$ 
in Euclidean space with no parabolic points oriented by a field of unit normals 
$n : M\to S^2$. 
Consider on $M$ the unique complex
structure compatible with the given orientation and the conformal structure defined
by the third fundamental form $\text{\rom{III}}=\text{d}n\cdot \text{d}n$.
(Here $\cdot$ denotes
the Euclidean inner product.)
Accordingly, the second fundamental form $\text{\rom{II}}$ decomposes into
bidegrees:
$$
\text{\rom{II}}=\text{\rom{II}}^{(2,0)}+\text{\rom{II}}^{(1,1)}+\text{\rom{II}}^{(0,2)}.
$$
$\text{\rom{II}}^{(2,0)}$ is a globally defined $(2,0)$ symmetric bilinear form on $M$
which plays the role of the usual Hopf differential for the pair of quadratic forms
$\text{\rom{III}}$, $\text{\rom{II}}$. We refer to $\text{\rom{II}}^{(2,0)}$ simply as
the Hopf differential.

\definition{Definition} A Riemann surface $M$ equipped with a 
holomorphic quadratic differential $q$ is called a {\it polarized surface}. An immersion 
$f : (M,q) \to \Bbb R^3$ is called $L$-{\it isothermic} if
$\text{\rom{II}}^{(2,0)}=\mu q$, for a real-valued smooth function $\mu : MÊ\to \Bbb R$.
\enddefinition 

\flushpar 
Near any point $p\in M$ where $q_{|p}\neq 0$ there exists local complex 
coordinate $z=x+iy : \Omega\subset M  \to \Bbb C$ such that $q_{\vert\Omega}=\text{d}z^2$
and 
$\text{\rom{III}}_{\vert\Omega}=e^{2u}\text{d}z\text{d}\overline{z}$. Then
$(x,y)$ are curvature line coordinates which are conformal for the third fundamental form
of $f$. $(x,y)$ will be called {\it conformal principal coordinates} for $f$ 
\footnote{The umbilic points of the immersion are precisely the zeros of the Hopf 
diffrential. The idea here is that not all umbilic points prohibit conformal principal
coordinates. If an umbilic is caused by a zero in the polarization, then it is a {\it bad}
umbilic, where no conformal principal coordinates may be found; instead, umbilic points 
caused by zeros in $\mu$ will not cause problems.}.

\medskip
\flushpar Examples of $L$-isothermic surfaces include 
surfaces of revolution, molding surfaces [BCG], surfaces with plane lines of curvature
in both systems [MN4], and the class of Weingarten surfaces on which $aH+bK=0$, for
constants $a,b$ with $a\neq 0$, where $H$ and $K\neq{0}$ denote the mean and Gaussian
curvatures,  respectively.
The last example follows as an application of Hopf's classical argument:
$H/K$ which is 
$H(\text{\rom{III}},\text{\rom{II}}):=(1/2)\text{tr}_{\text{\rom{III}}}\text{\rom{II}}$ is constant 
if and only if $\text{\rom{II}}^{(2,0)}$
is holomorphic [H].
\bigskip 
 
\subheading{3. The geometry of $L$-isothermic surfaces}  

\bigskip
\flushpar A pair of real quadratic forms $\text{\rom{III}}$ and $\text{\rom{II}}$ on
$M$ such that: $\text{\rom{III}}$ is positive definite, the intrinsic curvature 
$K(\text{\rom{III}})\equiv 1$, and such that $\text{\rom{II}}$ satisfies the Codazzi
equations with respect to the metric $\text{\rom{III}}$, 
defines, up to contact tansformations, an immersion 
$F=(f,n) : M\to \Bbb R^3\times S^2$ in the space of contact elements of $\Bbb R^3$  
satisfying the {\it contact condition} $\text{d}f\cdot{n}=0$. A map satisfying this
condition is called a Legendre immersion.
The {\it Euclidean projection} $f : M \to \Bbb R^3$ need not be an immersion nor even 
have constant rank.

\medskip 
\flushpar Let $\Lambda = \Bbb R^3\times S^2$ denote the space of contact 
elementsof $\Bbb R^3$.  
Now each contact element $(x,n)\in \Lambda$ corresponds to a null line
in Minkowski 4-space $\Bbb R^4_1$ with its Lorentz scalar product 
$\langle\ ,\ \rangle=(u^1)^2+ (u^2)^2+(u^3)^2-(u^4)^2$ via
$$
(x,n)\mapsto [x,n]:=
\left\{{^t(}x-tn,t), t\in \Bbb R\right\}\subset \Bbb R^4_1.
$$
Under this identification, 
the identity component $L =
\Bbb R^4 \rtimes SO_0(3,1)$ of the isometry group of $\Bbb R^4_1$ 
acts transitively on $\Lambda$ 
and preserves the contact condition. 
The action of $L$ permutes Legendre immersions: let $F=(f,n)$
be a Legendre immersion and, for $p\in M$, 
consider the null line $[f(p),n(p)]$. 
For each $A\in L$, $A\cdot[f(p),n(p)]\in \Lambda$ and
intersects $\Bbb R^3=\{{^t(}x,t)\in \Bbb R^4_1 :
t=0\}$ in the unique point $f'(p)$. If ${^t(}n',-1)$ denotes the
null direction of $A\cdot[f(p),n(p)]$, then $AF:= (f',n')$ defines another
Legendre immersion. Notice that the action does not preserve the Euclidean
fibration $\Lambda \to \Bbb R^3$.

\remark{Remark} A standard model for Laguerre geometry 
is obtained by identifying $\Bbb R^4_1$ with the set $\Cal S$ of 
oriented spheres (including point spheres) of $\Bbb R^3$ by
$$
\Sigma : \Cal S \to \Bbb R^4_1, \quad \sigma_r(p)\mapsto {^t(}p,r),
$$
where $\sigma_r(p)$ denotes a sphere with center $p$ and signed radius $r$.
Note that if $r>0$ (resp. $r<0$) then ${^t(}p,r)$ is the vertex of the 
backward (resp. forward) pointing light-cone which intersects $\Bbb R^3$ in 
exactly the sphere $\sigma_r(p)$. 
Two spheres $\sigma_r(p)$ and
$\sigma'_r(p)$ are in oriented contact if and only if 
${^t(}p-p',r-r')$ is a null
vector. Thus each contact element $(x,n)$ determines a null line of
2-spheres all of which are in oriented contact at $x$ with normal $n$. 
Following the classical terminology, a
{\it Laguerre transformation} is a contact transformation of $\Lambda$ 
induced by an
element in the group $L$. In terms of $\Bbb R^3$, a Laguerre
transformation takes oriented planes in $\Bbb R^3$ to
oriented planes, and oriented spheres to oriented spheres. 
In this context, 
two immersions $f,f' : M \to \Bbb R^3$ are
said to be {\it Laguerre equivalent} if there exists $A\in L$ such that 
their respective Legendrian 
lifts  $F, F'$ satisfy
$AF=F'$. (For more details about Laguerre geometry we refer to 
[Bl],[C],[MN1].) 
\endremark 

\medskip
\flushpar The theory of $L$-isothermic surfaces belongs in the Laguerre geometry: 
 
\proclaim{(Laguerre invariance)} Let $F=(f,n)$ be the Legendrian lift of an
$L$-iso\-ther\-mic surfaces $f$ and $A : \Lambda \to \Lambda$ be a Laguerre
contact transformation. Then the Euclidean projection of $AF=(f',n')$ is 
$L$-isothermic also.
\endproclaim

\demo{Proof} Let $\Cal L^{+}$ be the positive light-cone in $\Bbb R^4_1$. The
projective light-cone
$\Bbb P(\Cal L^{+})$ identifies with the conformal 2-sphere $S^2$ and the
projection $\Cal L^+ \to \Bbb P(\Cal L^+)$ 
is a principal $\Bbb R^+$--bundle which is trivial.
For each $A\in L$, $\text{d}A \in
SO_0(3,1)$, which preserves the light-cone $\Cal L^{+}$ and 
descends to an 
orientation-preserving conformal 
diffeomorphism $\tilde {A}$ of $S^2$.
This implies that $n'=\tilde {A}(n)$. Thus,
the conformal class of $\text{\rom{III}}$ and then the conformal
class of $\text{\rom{II}}$ (mod $\text{\rom{III}}$) are Laguerre invariant. 
\qed
\enddemo

\medskip
\remark{3.1 Notations: moving frames in Laguerre geometry}\endremark
\medskip

\noindent Here and in the following we shall consider $\Bbb R^4_1$ 
with linear coordinates
$x^1,\dots, x^4$ such that 
$\langle\ ,\ \rangle=-2x^1x^4 + (x^2)^2+(x^3)^2$; an 
orientation for which
$\text{d}x^1\wedge\dots\wedge\text{d}x^4\neq 0$; 
and a time-orientation given by 
$x^1+x^4>0$. By a {\it Laguerre frame} is meant a position vector
$a_0\in \Bbb R^4_1$ and an oriented basis $a_1,a_2,a_3,a_4$ of 
$\Bbb R^4_1$ such that
$(a_0;a_1,a_2,a_3,a_4)\in L$. Up to the choice of a 
reference frame, the manifold of
Laguerre frames may be identified with the group $L$. 
For any $A=(a_0,a)\in L$, $a_i$,
$i=1,2,3,4$, denote the column vectors of the matrix $a$.  
Geometrically, the null directions $a_1$, $a_4$
represent oriented planes which are in oriented contact with the 
oriented sphere
represented by $a_0$. By $[a_0,a_1], [a_0,a_4]$ we denote the lines 
in $\Bbb R^4_1$ through $a_0$ with null directions $a_1$ and $a_4$, 
respectively. 

\smallskip
\noindent The transitive action of $L$ on 
$\Lambda$ defines a principal $L_0$-bundle
$$
\pi_L : L \to\Lambda=L/L_0, 
 A \mapsto [a_0,a_1]. 
$$

\definition{Definition}  
A {\it Laguerre frame field} in $\Lambda$ is a local smooth 
section $A=(a_0,a)$ of $\pi_L$. 

\noindent A (local) {\it Laguerre framing along
a Legendre immersion} $F : M\to\Lambda$ is a smooth map $A : \Cal U \to L$
defined on an open subset $\Cal U\subset M$, such that
$\pi_LA = [a_0,a_1]= F$.
\enddefinition

\medskip
\remark{Remark (Principal frames)} Let $F : (M,q) \to \Lambda$ be a Legendrian immersion.
Then $F$ can be equipped with a Laguerre framing
$A=(a_0;a_1,a_2,a_3,a_4) : M\to L$
such that $F(p)=[a_0,a_1](p)$ and 
span$\{a_2,a_3\}_{|p}=\text{d}F_{|p}(T_pM)$, for any $p\in M$.
Let $\alpha = A^{-1}\text{d}A$ be its connection form. Then
$
\alpha^4_0=0,\  
\alpha^2_1\wedge\alpha^3_1\neq 0$.
Differentiating $\alpha^4_0=0$, it follows that the quadratic form  
$\alpha^2_0\alpha^2_1
+\alpha^3_0\alpha^3_1$ 
is symmetric and then diagonalizable. So we may assume that $\alpha^2_0\wedge\alpha^2_1=
\alpha^3_0\wedge\alpha^3_1=0$. We call such an $A$ a {\it principal framing} along $F$. 
An easy calculation using the Maurer--Cartan equation of $L$ yields
$$
\text{d}(\alpha^2_1 +i\alpha^3_1)=(\alpha^1_1 -*\alpha^3_2)\wedge(\alpha^2_1 +i\alpha^3_1),
$$
from which follows that {\it $F$ is $L$-isothermic if and only if
$\text{d}(\alpha^1_1 -\ast\alpha^3_2)=0$.}
\endremark

\bigskip

\subheading{4. Sphere congruences and the Bianchi--Darboux transformation}

\bigskip

\flushpar In classical surface theory, a sphere congruence is an immersion 
$S : M\to \Bbb R^4_1$ of a surface $M$ into the space $\Bbb R^4_1$ of all
oriented 2-spheres (including points) of $\Bbb R^3$.
A Legendre immersion $F=(f,n)$ is said to envelope the sphere congruence $S$ if
for each $p\in M$, the  sphere $\Sigma^{-1}(S(p))$ is in oriented contact 
with the Legendre surface at $F(p)$. 
If $S$ is a space-like immersion, i.e., the induced metric on $M$ is positive
definite, then there exist two enveloping surfaces [Bl]. It follows that 
the spheres
of the congruence are the common tangent spheres of the two envelopes.
Accordingly, there
results a mapping between the enveloping surfaces such that the sphere 
congruence
consists of the spheres passing through the points on the envelopes.

\definition{Definition}
A space-like sphere
congruence is called {\it Ribaucour} if the resulting mapping between the two
envelopes preserves curvature lines.
A Ribaucour sphere congruence is called {\it Darboux}
\footnote{In the context of M\"obius geometry a sphere congruence is an
immersion of a surface into the Lorentzian sphere $S^4_1$ interpreted as the 
space
of all oriented spheres (excluding points) and oriented planes in $\Bbb R^3$.
A sphere congruence is then called Ribaucour if the curvature lines on the two
envelopes correspond;   
a Ribaucour sphere congruence is Darboux if the correspondence between 
the two envelopes is conformal with respect to the first fundamental forms of
the envelopes, see for example [BHPP] and the literature therein.  } if the
resulting mapping is conformal with respect to  the third fundamental forms
induced on $M$ by the two envelopes $F$, $\hat{F}$. \enddefinition

\medskip
\flushpar Let $S : M \to \Bbb R^4_1$ be a space-like congruence 
enveloping $F$ and $\hat{F}$; we can 
consider an adapted  
Laguerre frame $A=(S;a_1,\dots,a_4) : M \to L$ such that
 $F=[S,a_1]$,  
$\hat{F} = [S,a_4]$, $a_2,a_3$ are tangential, and $a_1,a_4$ are normal
null directions.
The connection form $\alpha =(\alpha_0,\alpha')= A^{-1}\text{d}A$ of 
$A$ is then
$$
\left[\pmatrix
0\\
\alpha^2_0\\
\alpha^3_0\\
0\\
\endpmatrix,
\pmatrix
\format\c&\quad\c&\quad\c&\quad\c\\
 \alpha^1_1&\alpha^1_2&\alpha^1_3&0 \\
 \alpha^2_1&0&-\alpha^3_2&\alpha^1_2 \\
 \alpha^3_1&\alpha^3_2&0&\alpha^1_3 \\
 0&\alpha^2_1&\alpha^3_1&-\alpha^1_1
\endpmatrix\right]. 
$$
Observe that $\alpha^2_0\wedge\alpha^3_0\neq 0$ on $M$. 
The Maurer--Cartan equation $\text{d}\alpha +\alpha\wedge\alpha=0$ yields
the Ricci equations: 
$$
\align
&0=\alpha^2_0\wedge\alpha^1_2+\alpha^3_0\wedge\alpha^1_3, \tag1\\
&0=\alpha^2_0\wedge\alpha^2_1+\alpha^3_0\wedge\alpha^3_1, \tag2\\
&d\alpha^1_1+ \alpha^1_2\wedge\alpha^2_1+\alpha^1_3\wedge\alpha^3_1=0,\tag3
\endalign
$$
the Gauss equation:
$$
\text{d}\alpha^3_2=\alpha^1_2\wedge\alpha^3_1+\alpha^2_1\wedge\alpha^1_3,\tag4
$$
and the Codazzi equations:
$$
\align
&\text{d}\alpha^2_0=  \alpha^3_2\wedge\alpha^2_0,\quad
\text{d}\alpha^3_0=  -\alpha^3_2\wedge\alpha^2_0\tag5\\
&\text{d}\alpha^2_1= \alpha^1_1\wedge\alpha^2_1+\alpha^2_1\wedge\alpha^3_2,\quad
\text{d}\alpha^3_1= \alpha^1_1\wedge\alpha^3_1+\alpha^2_1\wedge\alpha^3_2\tag6\\
&\text{d}\alpha^1_2= \alpha^1_2\wedge\alpha^1_1+\alpha^3_2\wedge\alpha^1_3,\quad
\text{d}\alpha^1_3= \alpha^1_3\wedge\alpha^1_1-\alpha^3_2\wedge\alpha^1_2.\tag7
\endalign
$$

\medskip
We thus can state:

\proclaim{Theorem 1} Let $S : M\to \Bbb R^4_1$ be a flat space-like immersion
with flat normal bundle. Then $S$ is a Darboux sphere congruence
and both its enveloping surfaces -- which have opposite orientations -- are $L$-isothermic. 
\endproclaim

\demo{Proof} Assume $S$ induces a flat space-like metric on $M$,
i.e., $\text{d}\alpha^1_1=0$, and 
has flat normal bundle, i.e., $\text{d}\alpha^3_2=0$. 
According to (2), the second fundamental form $\alpha^2_0\alpha^2_1
+\alpha^3_0\alpha^3_1$ of $S$
in the normal direction $a_1$ is diagonalizable. So we may choose $A$ such that
$$
\alpha^2_0=h_2\alpha^2_1, \quad \alpha^3_0=h_3\alpha^3_1, \quad h_2-h_3\neq 0.\tag8
$$
In particular, $A$ becomes a principal framing along $F$. 
We can now write
$\alpha^1_2 =a_{11}\alpha^2_1+a_{12}\alpha^3_1$, 
$\alpha^1_3 =a_{21}\alpha^2_1+a_{22}\alpha^3_1$
for smooth functions $a_{ij}$. From   
equations (3) and (4) we obtain $a_{12}=a_{21}$,
$a_{11}=-a_{22}$. 
Further substituting into (1) and (2) yields
$$
a_{12}(h_{2}-h_{3})=0,
$$
and hence $a_{12}=0$. Therefore $\alpha^1_2 =a_{11}\alpha^2_1$ and   
$\alpha^1_3 =-a_{11}\alpha^3_1$. This implies, in particular, that 
the correspondence induced by $S$ on the two envelopes preserves 
curvature lines
and that $\langle\text{d}a_1,\text{d}a_1\rangle\varpropto
\langle\text{d}a_4,\text{d}a_4\rangle$, that is, $S$ is a Darboux sphere 
congruence.
Also, this tells us that the two envelopes have opposite orientations.

\smallskip
\noindent As for the isothermic property of the envelopes, by the 
Codazzi equations (7),
$$
\align
\text{d}a_{11}\wedge\alpha^2_1 &= -2a_{11}(\alpha^1_1\wedge\alpha^2_1+\alpha^3_2
\wedge\alpha^3_1)\\
\text{d}a_{11}\wedge\alpha^3_1 &= -2a_{11}(\alpha^1_1\wedge\alpha^3_1-\alpha^3_2
\wedge\alpha^2_1),
\endalign
$$
and from these
$2(\alpha^1_1 -*\alpha^3_2) = -\text{d}\log{|a_{11}|}$,
which is the condition for $F$ being $L$-isothermic according to the remark 
in the previous section. 

\smallskip 
\flushpar Concerning the second envelope $\hat{F}$,  
$B=(b_0;b_1,b_2,b_3,b_4):=$$(a_0;a_4,a_2,-a_3,a_1)$ defines a frame along 
$\hat{F}$. Its connection form is computed to be
$$
\beta = B^{-1}\text{d}B=
\left[\pmatrix
0\\
\alpha^2_0\\
-\alpha^3_0\\
0\\
\endpmatrix,
\pmatrix
\format\c&\quad\c&\quad\c&\quad\c\\
 -\alpha^1_1&\alpha^2_1&-\alpha^3_1&0 \\
 \alpha^1_2&0&\alpha^3_2&\alpha^2_1 \\
 -\alpha^1_3&-\alpha^3_2&0&-\alpha^3_1 \\
 0&\alpha^1_2&-\alpha^1_3&\alpha^1_1
\endpmatrix\right]. 
$$
$B$ is then a principal framing along $\hat{F}$, and the one form
$\beta^1_1 -*\beta^3_2$ is closed. So, also $\hat{F}$ is $L$-isothermic.  
\qed\enddemo

\medskip
\flushpar If $S$ is a Darboux congruence, then $S$ is either
as in Theorem 1, or
the normals $a_1$ and
$a_4$ of the two envelopes are M\"obius equivalent in $S^2$, that is 
$F$ and $\hat{F}$ are
Laguerre equivalent and have the same orientations. In the latter 
case the congruence $S$ belongs to a fixed sphere complex [Bl].

\medskip
Up to this degenerate situations, we can state:

\proclaim{Theorem 2} A map $S : M \to \Bbb R^4_1$ defines a
nondegenerate Darboux sphere congruence if and only if it is a flat space-like
immersion with flat normal bundle.
\endproclaim 
 
\definition{Definition} The $L$-isothermic surfaces enveloping a
non\-de\-gene\-rate Darboux congruence 
are called 
\it Bianchi--\-Darboux transforms \rm of each other. 
They form a Bianchi--Darboux pair. 
\enddefinition

\remark{Remark} Let $F=(f,n)$, $\hat{F}=(\hat{f},\hat{n})$ be a 
Bianchi--Darboux pair
enveloped by $S : M\to \Bbb R^4_1$. 
Next, let $\frak{so}(3,1) =\frak{k}\oplus\frak{m}$, 
$\frak{k}= \frak{so}(2)\times \frak{so}(1,1)$ be a symmetric 
decomposition of the Lie
algebra of $SO_0(3,1)$.  
Then, according to (3) and (4), the flatness of both the induced 
metric and the normal bundle of  $S$ is expressed by  
$$
\alpha'_{\frak{m}}\wedge\alpha'_{\frak{m}}=0,\tag9
$$
where $\alpha'_{\frak{m}}$ denotes the $\frak{m}$-part of $\alpha'$.
This condition expresses the fact that the map 
$(n,\hat{n}) : M\to N=S^2\times S^2\diagdown\Delta$ into the manifold
of ordered pairs of distinct points of $S^2$ is a curved flat, see [FP],[BHPP].
$N$ may be viewed as the Grassmannian $\tilde{G}_{1,1}(\Bbb R^4_1)$ of
oriented 2-planes of signature $(1,1)$ in $\Bbb R^4_1$ --- a pseudo-Riemannian
symmetric space with invariant metric of signature $(2,2)$. 
Note that the converse is not
true, that is, from a curved flat in $N$ we can only construct the normals
of a Bianchi--Darboux pair.

\endremark

\definition{Definition}
A map $A=(a_0,a) : M\to L$ such that $\alpha=(\alpha_0,\alpha')=A^{-1}dA$ 
satisfies (9) and $\alpha^1_0=\alpha^4_0=0$ is called a 
{\it curved flat framing} and $a_0 : M\to \Bbb R^4_1$ 
the associated sphere congruence. 
Note that, up to a gauge change, such a  frame can be chosen so that 
$\alpha^1_1=0$.
\enddefinition

The above discussion yields:

\proclaim{Corollary 1} The sphere congruence $a_0 : M\to \Bbb R^4_1$
of a curved flat framing $A=(a_0,a) : M\to L$ is a Darboux sphere
congruence enveloped by the two $L$-isothermic immersions
$F=[a_0,a_1]$ and $\hat{F}=[a_0,a_4]$.
Conversely, a (nondegenerate) Darboux congruence $S$ defines a curved 
flat framing $(S;a_1,a_2,a_3,a_4) : M \to L$,
where $\text{{\rm d}}S(TM)=\text{{\rm span}}\{a_2,a_3\}$ 
and $a_1$, $a_4$ generate the null subbundles of
the normal bundle of $S$. 
\endproclaim

\medskip

\flushpar In the next section we shall discuss the existence 
of Bianchi--Darboux transforms. 

\bigskip

\subheading{5. Integrability: construction of the Bianchi--Darboux transform}

\bigskip
 
\flushpar Let $F : M \to \Lambda$ be $L$-isothermic, and consider a 
principal frame
$A=$$(a_0;a_1,a_2,a_3,a_4)$ along $F$, $F=[a_0,a_1]$. Let $\alpha$ denote 
its connection 
form and $\rho$ be a smooth positive function such that
$$
2(\alpha^1_1 -*\alpha^3_2) = -\text{d}\log{\rho}.
$$
Define
$$
\alpha_\rho=(\alpha_0,\alpha'_\rho):=\left[\pmatrix
\alpha^1_0\\
\alpha^2_0\\
\alpha^3_0\\
0\\
\endpmatrix,
\pmatrix
\format\c&\quad\c&\quad\c&\quad\c\\
 \alpha^1_1&\alpha^1_2+\rho\alpha^2_1&\alpha^1_3-\rho\alpha^3_1&0 \\
 \alpha^2_1&0&-\alpha^3_2&\alpha^1_2+\rho\alpha^2_1 \\
 \alpha^3_1&\alpha^3_2&0&\alpha^1_3-\rho\alpha^3_1\\
 0&\alpha^2_1&\alpha^3_1&-\alpha^1_1
\endpmatrix\right]. 
$$
Note that $0=\text{d}\alpha_\rho + \alpha_\rho\wedge\alpha_\rho$.

\medskip
\flushpar Consider a solution 
$v={^t(}v^1,\dots,v^4) : M \to \Cal L^+\subset \Bbb R^4_1$
of the completely integrable linear system 
$$
\text{d}v=-{\alpha'_{\rho}}v.
$$
By a direct computation:

\proclaim{Lemma} {The one form 
$\quad{\gamma:=v^2 \alpha^2_0 + v^3\alpha^3_0 - v^4\alpha^1_0}\quad$ 
is closed.}  
\endproclaim

\flushpar Next put  
$s:={r\over{v_4}}$, where
\footnote{As it will be clear from the following discussion, the 
function $r$, which locally integrates $\gamma$, gives essentially the
signed radius of the Darboux sphere congruence we are going to 
construct (cf. Section 5.1).} $dr =\gamma$, and consider the gauged 
frame $\tilde A= Ag(s,v)$, 
where \footnote{Observe that the group of all elements $g(s,v)$, 
$s\in \Bbb R$,$v\in \Cal L^+$, is diffeomorphic to the structure
group of the fibration $\Cal P(F)\to M$ of principal frames along $F$.}   
$$
g(s,v)=\left[\pmatrix s\\
0\\
0\\
0\\
\endpmatrix,
\pmatrix
\format\c&\quad\c&\quad\c&\quad\c\\
 1/{v^4}&v^2/{v^4}&v^3/{v^4}&v^1\\
0&1&0&v^2 \\
 0&0&1&v^3 \\
 0&0&0&v^4  
\endpmatrix\right]. 
$$
Observe that $[\tilde{a_0},\tilde{a_1}]=[{a_0},{a_1}]=F$. 
The connection form 
$\tilde{A}^{-1}\text{d}\tilde{A}$ takes the form
$$
\!\!\left[\!\!\pmatrix
0\\
\alpha^2_0+s\alpha^2_1\\
\alpha^3_0+s\alpha^3_1\\
0\\
\endpmatrix,\!\!
\pmatrix
\format\c&\quad\c&\quad\c&\quad\c\\
 0&-\rho{v}^4\alpha^2_1&\rho{v}^4\alpha^3_1&0 \\
 {1\over{v^4}}\alpha^2_1&0&-\alpha^3_2+{1\over{v^4}}(v^3\alpha^2_1-v^2
\alpha^3_1)&-\rho{v}^4\alpha^2_1\\
 {1\over{v^4}}\alpha^3_1&\alpha^3_2-{1\over{v^4}}(v^3\alpha^2_1-v^2
\alpha^3_1)&0&
\rho{v}^4\alpha^3_1\\
 0&\alpha^2_1&\alpha^3_1&0
\endpmatrix\!\!\right]\!\!. 
$$
Thus $\tilde A$ is a curved flat framing and 
$\tilde{a}_0 : M\to \Bbb R^4_1$ is a Darboux
congruence which envelopes the Bianchi--Darboux pair $F$ and
$\hat{F}=[\tilde{a_0},\tilde{a_4}]$. 

We thus have proved:

\proclaim{Theorem 3}{Let $F : M \to \Lambda$ be an $L$-isothermic 
immersion and $A$ be any principal frame along $F$. 
Then any solution $v : M\to \Cal L^+$ of the linear system
$$
\text{{\rm{d}}}v=-{\alpha'_{\rho}}v
$$ 
defines an  $L$-isothermic immersion
$$
\hat{F}:=[a_0+sa_1,v^1a_1+v^2a_2+v^3a_3+v^4a_4] : M\to \Lambda
$$
which is a Bianchi-Darboux transform of $F$.}
\endproclaim

\medskip
\remark{Remark} The space of principal frames $\Cal P(F)$ along a 
Legendrian immersion $F$
can be viewed as a 6-dimensional integral submanifold of the
exterior differential system $\omega^4_0=0$, $\omega^2_0\wedge\omega^2_1=
\omega^3_0\wedge\omega^3_1=0$ on $L$ with independence condition 
$\omega^2_1\wedge\omega^3_1\neq 0$, where $\omega$ denotes
the Maurer--Cartan form of $L$. Proving the existence of Bianchi--Darboux
transforms for an $L$-isothermic $F$ can be reduced to checking the
Frobenius integrability condition for the Pfaffian system on 
$\Cal P(F)\times\Bbb R$ 
given by:
$$
\omega^1_0=0, \quad \omega^1_2-a\omega^2_1=0,\quad\omega^1_3+a\omega^3_1=0,
\quad \text{d}a +2a(\omega^1_1-*\omega^3_2)=0.
$$
\endremark

\medskip

\remark{5.1 The Bianchi-Darboux transform in Euclidean terms}\endremark

\medskip
\flushpar Let $f : (M,q) \to \Bbb R^3$ be an $L$-isothermic immersion 
with normal $n$ and
conformal principal coordinate $z = x+iy$.
If we identify $\Bbb E(3)$ with the subgroup of $L$
consisting of all $A\in L$ fixing the time-like vector $\epsilon_1+\epsilon_4$
($\epsilon_1,\dots,\epsilon _4$ the canonical basis 
of $\Bbb R^4_1$), 
the Euclidean framing $(f; t_1,t_2,t_3) : M \to \Bbb E(3)$ defined by 
$t_1 = n$, $t_2 =
{{{f}_x}\over{\|{f}_x\|}}$, $t_3= {{{f}_y}\over{\|{f}_y\|}}$ 
corresponds to the Laguerre
framing
 $$
E=(e_0;e)=
\left(
\pmatrix
{f^1}\over\sqrt{2}\\
f^2\\
f^3\\
-{f^1}\over\sqrt{2}
\endpmatrix;
\pmatrix
\format\c&\quad\c&\quad\c&\quad\c\\
 {1+t^1_1}\over2&{t^1_2}\over\sqrt{2}
&{t^1_3}\over\sqrt{2}&{1-t^1_1}\over2\\
{t^2_1}\over\sqrt{2}&t^2_2&t^2_3&-{t^2_1}\over\sqrt{2}\\
{t^3_1}\over\sqrt{2}&t^3_2&t^3_3&-{t^3_1}\over\sqrt{2}\\
{1-t^1_1}\over2&-{t^1_2}\over\sqrt{2}&-{t^1_3}\over\sqrt{2}
&{1+t^1_1}\over2\\
 \endpmatrix\right), 
$$
whose connection form 
can be written as
$$
\left[\pmatrix  
0\\
h_1\text{d}x\\
h_2\text{d}y\\
0\\
\endpmatrix,
\pmatrix
\format\c&\quad\c&\quad\c&\quad\c\\
 0&-{e^{u}\over{\sqrt{2}}}\text{d}x&-{e^{u}\over{\sqrt{2}}}\text{d}y&0 \\
 {e^{u}\over{\sqrt{2}}}\text{d}x&0&u_y\text{d}x -u_x\text{d}y&-{e^{u}
\over{\sqrt{2}}}\text{d}x\\
 {e^{u}\over{\sqrt{2}}}\text{d}y&-u_y\text{d}x +u_x\text{d}y&0&-{e^{u}
\over{\sqrt{2}}}\text{d}y \\%
 0&{e^{u}\over{\sqrt{2}}}\text{d}x&{e^{u}\over{\sqrt{2}}}\text{d}y&0%
\endpmatrix\right] 
$$  
for $u,h_1$, and $h_2$ smooth functions with $h_1h_2\neq 0$ at each point. 

\medskip
\flushpar In this setting $\rho = me^{-2u}$, for a constant  $m$, and 
the linear system becomes
$$
\left\{\matrix\format\l\\  
\text{d}v^1={{e^{u}-me^{-u}}\over{\sqrt{2}}}\text{d}xv^2+{{e^{u}+me^{-u}}
\over{\sqrt{2}}}\text{d}yv^3\\
\\
\text{d}v^2=-{e^{u}\over{\sqrt{2}}}\text{d}xv^1-(u_y\text{d}x -u_x
\text{d}y)v^3+{{e^{u}-me^{-u}}\over{\sqrt{2}}}\text{d}xv^4\\
\\
\text{d}v^3=-{e^{u}\over{\sqrt{2}}}\text{d}yv^1+(u_y\text{d}x -u_x
\text{d}y)v^2+{{e^{u}+me^{-u}}\over{\sqrt{2}}}\text{d}yv^4\\
\\
\text{d}v^4=-{e^{u}\over{\sqrt{2}}}\text{d}xv^1-{e^{u}\over{\sqrt{2}}}
\text{d}yv^2
\endmatrix
\right..\tag10
$$
Let ${^t(}v^1,v^2,v^3,v^4) : M\to \Cal L^+$ be a solution to (10) 
and let $r$ be a 
smooth function that locally integrates the closed 1-form 
$\gamma=v^2h_1\text{d}x + v^3h_2\text{d}y$;
we will refer to the functions $r, v^1, v^2, v^3, v^4$ as
{\it transforming functions}. We are in a position to state:

\proclaim{Theorem 4} Let $f : M\to \Bbb R^3$ be an $L$-isothermic immersion and
let $r,v^1,v^2,v^3,v^4$ be a set of transforming functions. 
Then
$$
\hat{f}=f+{{\sqrt{2}r}\over{v^1 +v^4}}{{{f}_x\times{f}_y}
\over{\|{{f}_x\times{f}_y}\|}}
- {{rv^2}\over{v^4(v^1 +v^4)}}
{{{f}_x}\over{\|{f}_x\|}}
-{{rv^3}\over{v^4(v^1 +v^4)}}{{{f}_y}\over{\|{f}_y\|}}.
$$
is a Bianchi--Darboux transform of $f$. All Bianchi--Darboux transforms of $f$
arise this way.
\endproclaim
\demo{Proof} 
Let $F=[e_0,e_1]$ be the Legendrian lift of $f$. According to Theorem 3, its
Bianchi--Darboux transform 
$\hat{F}=[e_0+{r\over{v^4}}e_1,v^1e_1+v^2e_2+v^3e_3+v^4e_4]$,
that is 
$$
\hat{F}=e_0+{r\over{v^4}}e_1+\mu(v^1e_1+v^2e_2+v^3e_3+v^4e_4),
$$ 
for a smooth function $\mu : M\to\Bbb R$. 
Now $\hat{F}$ takes values in $\Bbb E(3)$ if and only if
$\langle\hat{F},\epsilon_1+\epsilon_4\rangle=0$ if and only if 
$\mu={-r/{v^4(v^1+v^4)}}$. Substituting and using the above
realization of $\Bbb R^3$ in $\Bbb R^4_1$, we obtain the
required expression for $\hat{f}$. \qed\enddemo

\bigskip

\subheading{6. Superposition and permutability}

\bigskip
 
\flushpar Let $A^{(1)}, A^{(2)} : M \to L$ be two curved flat framings
having the same $L$-isothermic map $F =[a^{(1)}_0,a^{(1)}_1] =[a^{(2)}_0,a^{(2)}_1] : 
M\to \Lambda$
as first envelope, and with
second envelopes $F^{(1)}$ and $F^{(2)}$, respectively. Let $\alpha^{(1)}$, $\alpha^{(2)}$ 
be the corresponding connection forms.
 
\medskip
\flushpar The deformed forms associated with $\alpha^{(1)}$ are given by   
$$
\alpha^{(1)}_\lambda=
\left[\pmatrix
0\\
\alpha^2_0\\
\alpha^3_0\\
0\\
\endpmatrix,
\pmatrix
\format\c&\quad\c&\quad\c&\quad\c\\
 0&\lambda\alpha^1_2&\lambda\alpha^1_3&0 \\
 \alpha^2_1&0&-\alpha^3_2&\lambda\alpha^1_2 \\
 \alpha^3_1&\alpha^3_2&0&\lambda\alpha^1_3\\
 0&\alpha^2_1&\alpha^3_1&0
\endpmatrix\right]. 
$$
for some constant $\lambda \in \Bbb R$. Similarly for $\alpha^{(2)}$.  

\medskip

\flushpar  According to the discussion in Section 4,
 the framings $A^{(1)},A^{(2)}$ are related by a gauge
change 
$$ 
A^{(2)} = A^{(1)}g(s,v),\tag11
$$
where $v ={^t(}v^1,\dots,v^4) : M \to \Cal L^{+}$ is a solution to
the integrable linear system
$$
\text{d}v = -{\alpha^{(1)}_\lambda}' v,\tag12
$$
for some $\lambda$, and $s={r\over{v^4}}$ with 
$\text{d}r-v^2\alpha^2_0 -v^3\alpha^3_0 =0$.
Next, define  
the mapping $g^{-}_\lambda(s,v) : M \to L$ by 
$$
g^{-}_\lambda(s,v)=\left[\pmatrix
0\\
0\\
0\\
{{\lambda r}\over{v^1}}\\
\endpmatrix,
\pmatrix
\format\c&\quad\c&\quad\c&\quad\c\\
 \lambda^{-1}{v^1}&0&0&0\\
v^2&1&0&0\\
 v^3&0&1&0 \\
 \lambda v^4& {{\lambda v^2}\over{v^1}}&{{\lambda v^3}
\over{v^1}}&{{\lambda}\over{v^1}}
\endpmatrix\right]. 
$$

\flushpar Consider $\bar A^{(1)}:=A^{(1)}g^{-}_\lambda(s,v)$. 
The corresponding connection
form $\bar\alpha^{(1)}=(\theta,\eta)$ takes the form
$$
\align
&{^t\theta}= \left(
0,\ 
\alpha^2_0 +{{\lambda r}\over{v^1}}\alpha^1_2,\ 
\alpha^3_0 + {{\lambda r}\over{v^1}}\alpha^1_3,\ 
0\right),\\
&\eta=\pmatrix
\format\c&\quad\c&\quad\c&\quad\c\\
 0&{{\lambda}\over{v^1}}\alpha^1_2&{{\lambda}\over{v^1}}\alpha^1_3&0 \\
 { v^1}({1\over\lambda}-1)\alpha^2_1&0&-\alpha^3_2+{{\lambda v^3}\over{v^1}}\alpha^1_2
-{{\lambda v^2}\over{v^1}}\alpha^1_3
&{{\lambda}\over{v^1}}\alpha^1_2 \\
 { v^1}({1\over\lambda}-1)\alpha^3_1&
-\alpha^3_2+{{\lambda v^3}\over{v^1}}\alpha^1_2
-{{\lambda v^2}\over{v^1}}\alpha^1_3&0&{{\lambda}\over{v^1}}\alpha^1_3\\
 0&{ v^1}({1\over\lambda}-1)\alpha^2_1&{ v^1}({1\over\lambda}-1)\alpha^3_1&0
\endpmatrix.
\endalign 
$$

\flushpar We then have

\proclaim{Lemma} 
$\bar A^{(1)}:=A^{(1)}g^{-}_\lambda(s,v)$ is a curved flat framing
such that $[{\bar a^{(1)}}_0,{\bar a^{(1)}}_4]=F^{(1)}$.  
\endproclaim

 \flushpar Let $F'=[{\bar a^{(1)}}_0,{\bar a^{(1)}}_1]$ be the first 
envelope of the congruence 
${\bar a^{(1)}}_0$. 
We say that the $L$-isothermic map $F'$ is the \it superposition \rm of
$F^{(1)}$ and $F^{(2)}$ and write 
$$
F' = F^{(1)} \ast_{F}  F^{(2)}.
$$

\medskip
\proclaim{Theorem 5 (Permutability Theorem)} If an $L$-isothermic immersion 
$F$ has
two Bianchi--Darboux transforms $F^{(1)}$ and $F^{(2)}$, then there is another
$L$-isothermic immersion $F^{\ast}$ which is a Bianchi--Darboux transform of
$F^{(1)}$ and $F^{(2)}$ and is such that  
$$
F^{\ast}=F^{(1)} \ast_{F}  F^{(2)}=F^{(2)} \ast_{F}  F^{(1)}.
$$  
\endproclaim 

\demo{Proof} Write 
$$
A^{(1)}=A^{(2)}g(s,v)^{-1}=A^{(2)}g(\hat{s},\hat v),
$$
where
$$
\hat s = -s ,\quad \hat v = {^t(}v^1,-{{v^2}\over{v^4}},-{{v^3}
\over{v^4}},{1\over{v^4}}).\tag13
$$
By a direct calculation it is verified that, if $v : \Cal U\to \Cal L^{+}$
is a solution to (12), then $\hat v$ is a solution of 
$$
\text{d}\hat v = - {\alpha^{(2)}_\mu}' \hat v,
$$
where $\mu$ is given by
$$
\mu =\lambda(\mu - 1).\tag14
$$
Thus, $F^{(2)} \ast_{F}  F^{(1)}$ is the $L$-isothermic 
immersion represented by the first
envelope corresponding to $\bar A^{(2)}=A^{(2)}g^{-}_\mu(\hat{s},\hat{v})$.
It is now easily seen that, if
$\hat v$, $\hat{s}$, and $\mu$ are related to $v$, $s$, and 
$\lambda$ as in (13) and (14), then
$[{\bar a^{(2)}}_0, {\bar a^{(2)}}_1]=[{\bar a^{(1)}}_0, {\bar a^{(1)}}_1]$. 
The situation is visualized in the following diagram: 

$$
\CD
F            @>  A^{(2)} >>                F^{(2)}\\
@VV  A^{(1)} V                              @VV \bar A^{(2)} V\\
F^{(1)}            @> \bar A^{(1)}>>         F^{(1)} \ast_{F}  F^{(2)}
\endCD.
$$
\qed\enddemo

\enddocument